\documentclass[a4paper,12pt]{amsart}

\usepackage{amsfonts}
\usepackage{amssymb}
\usepackage{amsbsy}

\begin{document}

\title[Hopf-Galois Structures]{Nilpotent and Abelian Hopf-Galois
  Structures on Field Extensions} 
 
\author{Nigel P.~Byott}

\address{College of Engineering,
  Mathematics and Physical Sciences, University of Exeter, Exeter 
EX4 4QF U.K.}  \email{N.P.Byott@exeter.ac.uk}

\date{\today}

\newcommand{\Z}{{\mathbb{Z}}}
\newcommand{\Q}{{\mathbb{Q}}}
\newcommand{\F}{{\mathbb{F}}}

\newcommand{\gp}{{\mathfrak{p}}}
\newcommand{\gP}{{\mathfrak{P}}}

\newcommand{\sP}{{\mathcal{P}}}
\newcommand{\sW}{{\mathcal{W}}}

\newcommand{\Gal}{\mathrm{Gal}}

\newcommand{\lra}{\longrightarrow}
\newcommand{\bs}{\backslash}

\newtheorem{theorem}{\sc Theorem}
\newtheorem{proposition}[subsection]{\sc Proposition}
\newtheorem{corollary}[subsection]{\sc Corollary}
\newtheorem{lemma}[subsection]{\sc Lemma}
\newtheorem{example}[subsection]{\sc Example}
\newtheorem{remark}[subsection]{\sc Remark}

\newenvironment{pf}{{\em Proof. }}{\qed \newline}
\newenvironment{pf-ab-thm}{{\noindent \em Proof of Theorem \ref{ab-thm}. }}{\qed \newline}
\newenvironment{pf-nilp-thm}{{\noindent \em Proof of Theorem
    \ref{nilp-thm}. }}{\qed \newline}
\newenvironment{pf-cyc-thm}{{\noindent \em Proof of Theorem
    \ref{cyc-thm}. }}{\qed \newline}
\newenvironment{pf-all-Gamma-thm}{{\noindent \em Proof of Theorem \ref{all-Gamma-thm}. }}{\qed \newline}

\newcommand{\Aut}{\mathrm{Aut}}
\newcommand{\GL}{\mathrm{GL}}
\newcommand{\Hol}{\mathrm{Hol}}
\newcommand{\Perm}{\mathrm{Perm}}

\newcommand{\eab}{e_\mathrm{ab}}
\newcommand{\enil}{e_\mathrm{nil}}

\newcommand{\N}{\mathbb{N}}

\newcommand{\onto}{\twoheadrightarrow}

\maketitle{}

 

\begin{abstract}
Let $L/K$ be a finite Galois extension of fields with group $\Gamma$.
When $\Gamma$ is nilpotent, we show that the problem of enumerating
all nilpotent Hopf-Galois structures on $L/K$ can be reduced to the
corresponding problem for the Sylow subgroups of $\Gamma$. We use this
to enumerate all nilpotent (resp.~abelian) Hopf-Galois structures on a
cyclic extension of arbitrary finite degree. When $\Gamma$
is abelian, we give conditions under which every abelian
Hopf-Galois structure on $L/K$ has type $\Gamma$. 
We also give a criterion on $n$ such that \emph{every} Hopf-Galois structure
on a cyclic extension of degree $n$ has cyclic type. 
\end{abstract}

\section{Introduction and statement of results}

Let $\Gamma$ be a finite group and let $L/K$ be a finite extension of
fields with $\Gal(L/K) \cong \Gamma$ (for brevity, we say: $L$ is a
$\Gamma$-extension of $K$). Then $L$ is a module over the group algebra
$K[\Gamma]$, and $K[\Gamma]$ carries the structure of a $K$-Hopf algebra. This
makes $L$ into a $K[\Gamma]$-Hopf-Galois extension of $K$. There may be
other $K$-Hopf algebras $H$ which act on $L$ so that $L$ is an
$H$-Hopf-Galois extension. Such Hopf-Galois structures were investigated by
Greither and Pareigis \cite{GP}, who showed how the determination of
all Hopf-Galois structures on a given separable field extension $L/K$
could be reduced to a question in group theory. In particular, any
Hopf algebra $H$ which gives a Hopf-Galois structure on $L$ has the
property that $L \otimes_K H = L[G]$ as $L$-Hopf algebras, where $G$
is some regular group of permutations of $\Gamma$. Thus $G$ and
$\Gamma$ have the same order, but in general they need not be
isomorphic. We will refer to the isomorphism class of $G$ as the {\em type} of
the Hopf-Galois structure, and will say that the Hopf-Galois
structure is {\em abelian} (resp.~{\em nilpotent}) if $G$ is abelian
(resp.~nilpotent). 

For some groups $\Gamma$ it is known that every Hopf-Galois structure
on a $\Gamma$-extension must have type
$\Gamma$. This holds for cyclic groups of order $p^n$ with $p>2$
prime and $n \geq 1$ \cite{K}, for elementary abelian groups of
order $p^2$ with $p>2$ \cite{unique}, for cyclic groups of order $n$
with $(n, \varphi(n))=1$ (where $\varphi$ is Euler's totient function)
\cite{unique}, and for non-abelian simple groups \cite{simple}. On the
other hand, there are many groups $\Gamma$ for which there are 
Hopf-Galois structures whose type is different from $\Gamma$, the smallest
cases being the two groups of order 4 \cite{unique}. Indeed, if
$\Gamma$ is abelian then there may be Hopf-Galois structures which are
not abelian, or even nilpotent. For example, if $\Gamma$ is cyclic of order
$pq$, where $p$, $q$ are primes such that $q|(p-1)$, then $L/K$ admits
$2(q-1)$ Hopf-Galois structures which are not nilpotent, in addition
to the unique (classical) one of type $\Gamma$ \cite{pq}.
This phenomenon was investigated in some detail in \cite{NYJM}, where
it was shown that any abelian extension $L/K$ of even degree $n>4$
admits a non-abelian Hopf-Galois structure, and that the same holds
for many abelian groups of odd order. On the other hand, some new
groups $\Gamma$ were given in \cite{NYJM} for which all
Hopf-Galois structures are of type $\Gamma$ (cf.~Remark \ref{not-best}
below).   
 
In this paper, we supplement the results of \cite{NYJM} by considering
the situation where $\Gamma$ and $G$ are both abelian or, more
generally, both nilpotent. We will show that the enumeration of such
Hopf-Galois structures can be reduced to the case of groups of prime
power order.

Let $e(\Gamma,G)$ denote the number of Hopf-Galois structures of type
$G$ on a $\Gamma$-extension $L/K$. Thus the total number of 
Hopf-Galois structures on $L/K$ is given by
$$ e(\Gamma) = \sum_{G} e(\Gamma,G), $$
where the sum is over all isomorphism classes of groups $G$ of order
$|\Gamma|$. We also write
$$ \eab(\Gamma) = \sum_{G\ \small{\mathrm{abelian}}} e(\Gamma,G), \qquad 
   \enil(\Gamma) = \sum_{G\ \small{\mathrm{nilpotent}}}
     e(\Gamma,G), $$
where the sum is over all isomorphism types of abelian
(resp.~nilpotent) groups $G$ of order $|\Gamma|$. Thus $\eab(\Gamma)$
(resp.~$\enil(\Gamma)$) is the number of abelian (resp.~nilpotent)
  Hopf-Galois structures on $L/K$. Recall that a finite group $\Delta$
  is nilpotent if it is the 
direct product of its Sylow subgroups \cite[(5.2.4)]{Rob}. In
particular, if $\Delta$ is abelian, or if $\Delta$ is a $p$-group for some prime
number $p$, then $\Delta$ is nilpotent.

Let $n$ be the degree of the extension $L/K$. We write the prime
factorisation of $n$ as  
$$   n = \prod_{p | n} p^{v_p}, $$
where the product is over the distinct prime factors $p$ of $n$. If
$\Gamma$ is nilpotent, we can correspondingly write $\Gamma$ as a direct
product of groups
\begin{equation} \label{nilp-decomp}
  \Gamma = \prod_{p | n} \Gamma_p, 
\end{equation}
where $\Gamma_p$ is the (unique) Sylow $p$-subgroup of $\Gamma$ and
has order $p^{v_p}$. By Galois theory, we can then decompose $L$ as 
$$ L = \bigotimes_{p | n } L_p, $$ (tensor product over $K$) where
$L_p$ is a $\Gamma_p$-extension of $K$. If, for each $p$, we
take a Hopf-Galois structure on 
$L_p/K$, say of type $G_p$ and with corresponding $K$-Hopf algebra
$H_p$, then the Hopf algebra $H=\bigotimes_{p | n} H_p$ acts in the
obvious way on $L$, 
giving $L/K$ a Hopf-Galois structure of type $G= \prod_{p | n}
G_p$. This Hopf-Galois structure is necessarily nilpotent,
and is abelian if and only if each $G_p$ is abelian.

We will see that if $\Gamma$ is nilpotent then {\em every} nilpotent 
Hopf-Galois structure on $L/K$ arises in this way. This is the key
observation in the proof of our first main result:

\begin{theorem} \label{nilp-thm}
Let $\Gamma$ be a nilpotent group of order $n$. Then for each
nilpotent group $G$ of order $n$ we have $e(\Gamma,G)= \prod_{p | n}
e(\Gamma_p,G_p)$.
\end{theorem}

Taking the sum over all isomorphism types of nilpotent (resp.~abelian) 
groups $G$ of order $n$, we immediately obtain:

\begin{corollary} \label{nilp-ab}
For a finite nilpotent group $\Gamma$, we have 
$$ \enil(\Gamma)=\prod_{p | n} e(\Gamma_p) \mbox{ and } 
  \eab(\Gamma)=\prod_{p | n} \eab(\Gamma_p).  $$
\end{corollary}

As an application of Theorem \ref{nilp-thm}, we will determine
the number of nilpotent (resp.~abelian) Hopf-Galois structures 
on a cyclic extension of arbitrary finite degree. Before stating the
result, we fix some notation. For $m \geq 1$, let $C_m$ 
denote the cyclic group of order $m$, and, for $v \geq 3$, let $D_{2^v}$
  (resp.~$Q_{2^v}$)  denote the dihedral (resp.~generalized
  quaternion) group of order $2^v$. Also, for $n \geq 1$, let 
$r(n)$ be the radical of $n$:
$$ r(n) = \prod_{p | n} p. $$
\begin{theorem} \label{cyc-thm}
Let $\Gamma$ be a cyclic group of order $n$.
\begin{itemize}
\item[(i)]
If $n$ is not divisible by $4$, then
$$ \enil(\Gamma)=\eab(\Gamma)=e(\Gamma,\Gamma)= \frac{n}{r(n)}. $$ 
Thus every nilpotent Hopf-Galois structure on a cyclic extension of
degree $n$ is cyclic, and hence abelian.
\item[(ii)] If $n \equiv 4 \pmod{8}$, then again
$$ \enil(\Gamma)=\eab(\Gamma)= \frac{n}{r(n)}, $$ 
but
$$ e(\Gamma,\Gamma)=e(\Gamma, C_2 \times C_{n/2}) =
\frac{n}{2r(n)}. $$
Thus every nilpotent Hopf-Galois structure on a cyclic extension of
degree $n$ is abelian, but only half of them are cyclic.
\item[(iii)] If $n$ is divisible by 8, so $n=2^v n'$ with $v \geq 3$
  and $n'$ odd, then 
$$ \enil(\Gamma)= \frac{3n}{2r(n)} \mbox{ and } 
  \eab(\Gamma)=e(\Gamma,\Gamma)= \frac{n}{2r(n)},  $$ 
with 
$$ e(\Gamma,D_{2^v} \times C_{n'}) = e(\Gamma,Q_{2^v} \times C_{n'}) =
\frac{n}{2r(n)},  $$   
Thus every abelian Hopf-Galois structure on a cyclic extension of
degree $n$ is cyclic, although there are also Hopf-Galois structures which
are nilpotent but not abelian.
\end{itemize}
\end{theorem}

For a finite abelian $p$-group $\Gamma$, Featherstonhaugh, Caranti and
Childs \cite{FCC} have given conditions under which every abelian 
Hopf-Galois structure on a $\Gamma$-extension must have type
$\Gamma$. Combining this with Theorem \ref{nilp-thm}, we will obtain
the following result in the abelian case.

\begin{theorem} \label{ab-thm}
Let $\Gamma$ be a finite group of order $n = \prod_p p^{v_p}$, and
suppose that, for each prime factor $p$ of $n$, either $v_p <
p-1$ or $p \leq3$, $v_p < p$. Then every abelian Hopf-Galois structure
on a $\Gamma$-extension has type $\Gamma=\Gal(L/K)$. Equivalently,
$\eab(\Gamma)=e(\Gamma,\Gamma)$.
\end{theorem}

Combining Theorems \ref{cyc-thm} and \ref{ab-thm} with a result of
L.~E.~Dickson \cite{Dickson} dating from 1905, we obtain some new
cyclic groups $\Gamma$ for which {\em every} Hopf-Galois structure has
type $\Gamma$:

\begin{theorem} \label{all-Gamma-thm}
Suppose that $n= \prod_p p^{v_p}$ satisfies the following conditions:
\begin{itemize}
\item[(i)] $v_p \leq 2$ for all primes $p$ dividing $n$;
\item[(ii)] $p \nmid (q^{v_q}-1)$ for all primes $p$, $q$ dividing
  $n$;
\item[(iii)] $ 4 \nmid n$.
\end{itemize}
Then a cyclic extension of degree $n$ admits precisely
$n/r(n)$ Hopf-Galois structures, all of which are of cyclic type.
\end{theorem}

\medskip

\noindent {\sc Acknowledgment:} The author thanks Lindsay Childs and Tim Kohl
for email correspondence about this work, which led to a
simplification of some of the arguments. 

\section{Nilpotent Hopf-Galois Structures}

In this section we prove Theorem \ref{nilp-thm}.

We first recall the method of counting Hopf-Galois
structures on a $\Gamma$-extension for an arbitrary finite group $\Gamma$.
It was shown in \cite{GP} that these Hopf-Galois structures
correspond to regular permutation groups on $\Gamma$ which are
normalized by the group $\lambda(\Gamma)$ of left multiplications by
elements of $\Gamma$. (Recall that a permutation group $H$ on a set $X$
is regular if, given $x$, $y \in X$, there is a unique $h \in H$ with
$hx=y$.) Thus finding all Hopf-Galois structures with a given type $G$
amounts to finding all regular subgroups in the group
$\Perm(\Gamma)$ of permutations of $\Gamma$ which are isomorphic
to $G$ and are normalized by $\lambda(\Gamma)$. It was shown in
\cite{unique} that this problem can be reframed as a calculation
inside $\Hol(G) = \rho(G) \cdot \Aut(G)$, the holomorph of $G$,
which is usually a much smaller group than $\Perm(\Gamma)$. Here $\rho\;
\colon G \lra \Perm(G)$ is the right regular representation
$\rho(g)(x)=xg^{-1}$ for $g$, $x \in G$. As further reformulated by
Childs (see e.g.~\cite[\S7]{Ch00}), this gives the following method
of counting Hopf-Galois structures. A homomorphism  
$\beta \colon \Gamma \lra \Hol(G)$ will be called a regular embedding if it is
injective and its image is a regular group of permutations on
$G$. Two such embeddings will be called equivalent if they are conjugate by
an element of $\Aut(G)$. Then the number $e(\Gamma,G)$ of Hopf-Galois
structures of type $G$ on a $\Gamma$-extension
is the number of equivalence classes of regular
embeddings of $\Gamma$ into $\Hol(G)$. 

We will need the following general result.

\begin{proposition} \label{cent-prop}
Let $N$ be a regular subgroup of $\Hol(G)$. Then the centralizer of $N$
in $\Hol(G)$ has order dividing $|G|$. 
\end{proposition}
\begin{pf}
We can regard $\Hol(G)$ as a subgroup of the group $B=\Perm(G)$ of all
permutations of $G$. By \cite[Lemma 2.4.2]{GP}, the centralizer of $N$
in $B$ is canonically identified with the opposite group of $N$, so in
particular has order $|N|=|G|$. The centralizer of $N$ in $\Hol(G)$ is a
subgroup of this, so has order dividing $|G|$. 
\end{pf}

If $G$ is a nilpotent group, its Sylow subgroups $G_p$ are
characteristic subgroups. We therefore have direct product decompositions
\begin{equation} \label{dec-Aut}
   \Aut(G) = \prod_{p | n} \Aut(G_p), 
\end{equation}
and hence
\begin{equation} \label{dec-Hol} 
   \Hol(G) = \prod_{p | n} \Hol(G_p). 
\end{equation}

Now suppose that $\Gamma$ and $G$ are nilpotent groups of order $n$,
and that we are given a homomorphism $\beta_p \colon \Gamma_p \lra
\Hol(G_p)$ for each $p|n$. Using (\ref{nilp-decomp}) and (\ref{dec-Hol}), we 
can define a homomorphism
\begin{equation} \label{dp-beta}
  \beta = \left(\prod_{p|n} \beta_p\right) \colon \Gamma \lra \Hol(G).
\end{equation}
It is clear that if each $\beta_p$ is a regular embedding then so is
$\beta$. This construction corresponds to taking tensor products of
Hopf-Galois structures on field extensions of prime-power degrees, as
described in \S1.  

Not every homomorphism $\beta \colon \Gamma \lra \Hol(G)$ arises as
such a  product. For any primes $p$, $q$ dividing $n$, let
$\iota_p \colon \Gamma_p \lra \Gamma$ be the inclusion induced by the
direct product decomposition (\ref{nilp-decomp}) of $\Gamma$, and let
$\pi_q \colon \Hol(G) \lra \Hol(G_q)$ be the projection induced by
(\ref{dec-Hol}). Given a homomorphism $\beta \colon \Gamma \lra
\Hol(G)$, let $\beta_{pq}$ be the composite homomorphism $\beta_{pq}=\pi_q \circ
\beta \circ \iota_p \colon \Gamma_p \lra \Hol(G_q)$. Then $\beta$ is
determined by its matrix of components $(\beta_{pq})$. For each $q$,
the images of the $\beta_{pq}$ must centralize each other in
$\Hol(G_q)$, since the $\Gamma_p$ centralize each other in $\Gamma$.  
Conversely, a matrix of homomorphisms
$(\beta_{pq})$, $\beta_{pq} \colon \Gamma_p \lra \Hol(G_q)$, 
determines a homomorphism $\beta \colon \Gamma \lra \Hol(G)$, provided
only that, for each $q$, the images of the $\beta_{pq}$ centralize
each other in $\Hol(G_q)$. 

We can determine from the matrix $(\beta_{pq})$ whether $\beta$
is a regular embedding:

\begin{lemma} \label{diag}
Let $\Gamma$ and $G$ be nilpotent, and let $\beta \colon \Gamma \lra
G$ correspond to the matrix of homomorphisms $(\beta_{pq})$ as above.
Then $\beta$ is a regular embedding if and only if
$\beta_{pp} \colon \Gamma_p \lra \Hol(G_p)$ is a regular embedding for each $p$.
\end{lemma}
\begin{pf}
First observe that $\beta_{pp}(\Gamma_p)$ is the unique Sylow
$p$-subgroup in the subgroup $\pi_p \circ \beta(\Gamma)$ of $\Hol(G_p)$,
and hence is normal in $\pi_p \circ \beta(\Gamma)$. 

If $\beta$ is regular then $\pi_p \circ \beta(\Gamma)$ is transitive
on $G_p$. Then, by Proposition \ref{perm-prop} below, the number of
orbits of $\beta_{pp}(\Gamma_p)$ on $G_p$ divides both $|G_p|=p^{v_p}$
and $|\pi_p \circ \beta(\Gamma)/\beta_{pp}(\Gamma)|$ (which is coprime
to $p$).  Thus $\beta_{pp}$ is transitive, and hence regular, on
$G_p$.

Conversely, suppose that each $\beta_{pp}$ is a regular embedding. We write
$e_G$ for the identity element of $G$.  Consider the subsets
$X=\beta(\Gamma) e_G$ and $Y=\beta(\Gamma_p) e_G$ of $G$. Clearly $|Y|\leq
|\Gamma_p|$, and the regularity of $\beta_{pp}$ ensures that $|Y| \geq
|G_p|=|\Gamma_p|$. Hence $|Y|=|\Gamma_p|$. As $\beta(\Gamma_p)$ is
normal in $\beta(\Gamma)$, Proposition \ref{perm-prop} shows that all
orbits of $\beta(\Gamma_p)$ on $X$ have the same size. One such orbit
is $Y$, so $|X|$ is divisible by $|\Gamma_p|$. This holds for all $p$,
so $X=G$ and $\beta$ is a regular embedding.
\end{pf}
 
In the above proof, we used the following simple fact about
permutation groups:

\begin{proposition} \label{perm-prop}
Let $H$ be a finite group acting transitively on a set $X$, and let
$N$ be a normal subgroup of $H$.  Then the orbits of $N$ on $X$ all
have the same size, and the number of these orbits divides both 
$|X|$ and $|H/N|$.
\end{proposition}
\begin{pf}
Let $N$ have $m$ orbits on $X$, and let $Nx$ and $Ny$ be two such
orbits. Then $y=hx$ for some $h\in H$, and $Ny=Nhx=hNx$.  This shows
that the quotient group $H/N$ acts on the set $\{Nx\}$ of orbits of
$N$, and that this action is transitive. It follows firstly that these
orbits have the same size, so that $m$ divides $|X|$, and secondly that
$m$ divides $|H/N|$.
\end{pf}

\begin{pf-nilp-thm}
Let $\beta \colon \Gamma \lra \Hol(G)$ be a regular embedding, and let
$(\beta_{pq})$ be the corresponding matrix of homomorphisms. By Lemma
\ref{diag}, each $\beta_{pp}$ is a regular embedding of $\Gamma_p$
into $\Hol(G_p)$. For $p \neq q$,
the image of the homomorphism $\beta_{pq} \colon \Gamma_p \lra
\Hol(G_q)$ must centralize the regular subgroup $\beta_{qq}(\Gamma_q)$
of $\Hol(G_q)$, and so must be a $q$-group by Proposition
\ref{cent-prop}. But $\beta_{pq}(\Gamma_p)$ is a $p$-group since
$\Gamma_p$ is. Thus $\beta_{pq}$ is the trivial homomorphism
whenever $p \neq q$.  This means that the matrix $(\beta_{pq})$ is
``diagonal'' and $\beta$ is just the product $\beta= (\prod_p
\beta_{pp})$ as in (\ref{dp-beta}). Conversely, given a regular
embedding $\beta_p \colon \Gamma_p \lra \Hol(G_p)$ for each $p$, the
homomorphism $(\prod_p \beta_p) \colon \Gamma \lra G$ is a regular
embedding. It is immediate that these two constructions are mutually
inverse.

We have just established a bijection between regular embeddings 
$\beta \colon \Gamma \lra \Hol(G)$ and families of regular embeddings
$\beta_p \colon \Gamma_p \lra \Hol(G_p)$ for each $p | n$. It follows 
from (\ref{dec-Aut}) that two regular embeddings $\beta$, $\beta'$
are conjugate by an element of $\Aut(G)$ if and only if, for each $p$, their
components $\beta_p$, $\beta'_p$ are conjugate by an element of
$\Aut(G_p)$. Hence the equivalence classes of
regular embeddings $\beta \colon \Gamma \lra \Hol(G)$ correspond
bijectively to families of equivalence classes of regular embeddings
$\beta_p \colon \Gamma_p \lra \Hol(G_p)$. This shows that 
$e(\Gamma,G)=\prod_p e(\Gamma_p,G_p)$.
\end{pf-nilp-thm}

\section{Hopf-Galois structures on cyclic extensions} 

For cyclic extensions whose degree is a power of a prime $p$, all the 
Hopf-Galois structures are already known. We recall the results.  

\begin{lemma} \label{kohl-etc}

\begin{itemize}
\item[(i)] For $n=p^v$ with $p>2$ and $v\geq 1$, we have
  $e(C_n)=e(C_n,C_n)=p^{v-1}$. 
\item[(ii)] For $n=2$, we have $e(C_2)=e(C_2,C_2)=1$; for $n=4$, we
  have $e(C_4)=2$ with $e(C_4,C_4)=e(C_4,C_2\times C_2)=1$. 
\item[(iii)] For $n=2^v$ with $v \geq 3$, we have $e(C_n)=3 \cdot
  2^{v-2}$ with $e(C_n,C_n)=e(C_n,D_n)=e(C_n,Q_n)=2^{v-2}$. 
\end{itemize}
Thus, for a prime power $n=p^v$, we have $e(C_n)=n/r(n)$ except in the
case $p=2$, $v \geq 3$, when $e(C_n)=3n/(2r(n))$.  
\end{lemma}
\begin{pf}
(i) is equivalent to Kohl's result \cite{K} that, for an odd prime $p$, a cyclic Galois
  extension of degree $p^r$ admits $p^{r-1}$ Hopf-Galois
  structures, all of cyclic type. Similarly, (ii) follows from
  \cite{unique} and (iii) from \cite{2power}.
\end{pf}

Theorem \ref{cyc-thm} follows directly from Lemma
\ref{kohl-etc} and Theorem \ref{nilp-thm}.
 
\section{Abelian Hopf-Galois Structures}

In this section, we prove Theorems \ref{ab-thm} and
\ref{all-Gamma-thm}.

From \cite[Theorem 1]{FCC} we have the following result:

\begin{lemma} \label{FCC-lemma}
Let $\Gamma$ be an abelian $p$-group of $p$-rank $m$, with $p>m+1$. Then
$\eab(\Gamma)=e(\Gamma,\Gamma)$.  
\end{lemma}

\begin{pf-ab-thm}  
 Let $G$ be an abelian group of order $n$, and let $\Gamma_p$,
 $G_p$ be the Sylow $p$-subgroups of $\Gamma$, $G$ as usual. If
 $v_p<p-1$ then certainly $p>m+1$ where $m$ is the $p$-rank of $G_p$,
 so, by Lemma \ref{FCC-lemma}, $e(\Gamma_p,G_p)=0$ unless
 $G_p=\Gamma_p$.  If $p=3$ and $v_3=2$ then either $\Gamma_3 =
 C_9$, when by Lemma \ref{kohl-etc}(i) we have $e(\Gamma_3,G_3)=0$
 unless $G_3=\Gamma_3$, or $\Gamma_3 = C_3 \times C_3$, when the
 same conclusion holds by \cite{unique}. If $p=2$ and $v_2=1$ then
 $\Gamma_2 = C_2$ and $G_2 = C_2$. Thus the hypotheses of
 Theorem \ref{ab-thm} ensure that
 $\eab(\Gamma_p)=e(\Gamma_p,\Gamma_p)$ for all $p$. By Corollary
 \ref{nilp-ab} we then have
$$ \eab(\Gamma) = \prod_{p | n} e(\Gamma_p, \Gamma_p) =
e(\Gamma,\Gamma), $$
and every abelian Hopf-Galois structure on $L/K$ has type
$\Gamma$.
\end{pf-ab-thm}

To prove Theorem \ref{all-Gamma-thm}, we need the following old result of
L.~E.~Dickson \cite{Dickson} (see also \cite[\S5.5, Exercise 24,
  p.~189]{DF}):

\begin{lemma} \label{dickson}
Let $n$ have prime factorisation $\prod_p p^{v_p}$. Then every group
of order $n$ is abelian if and only if $v_p \leq 2$ for each prime $p$
dividing $n$, and $p \nmid (q^{v_q}-1)$ for all primes $p$, $q$
dividing $n$.
\end{lemma}

\begin{pf-all-Gamma-thm}
Let $\Gamma$ be a cyclic group of order $n$. 
The conditions of Theorem \ref{all-Gamma-thm} imply those of Theorem
\ref{ab-thm}, so that every abelian Hopf-Galois structure on a 
$\Gamma$-extension has cyclic type. On the other hand,
the hypotheses of Lemma \ref{dickson} are also satisfied. Thus every
group of order $n$ is abelian, and therefore every Hopf-Galois structure is
abelian. It follows that all the Hopf-Galois structures are cyclic. By
Theorem \ref{cyc-thm}(i), the number of Hopf-Galois structures is
therefore $n/r(n)$.  
\end{pf-all-Gamma-thm}

\begin{remark} \label{not-best}
In Theorem \ref{all-Gamma-thm}, there are no non-abelian Hopf-Galois
structures for the rather trivial reason that there are no non-abelian
groups of the appropriate order. This result is certainly not best
possible, since if $n=p^2 q^2$ for primes $2<q<p$ with $(q,p+1)>1$
(e.g.~$q=3$, $p=11$), or if $n=p^3 q$ for distinct primes $p$, $q$ with
$(p,q-1)=(q,p^2-1)=1$ but $(q,p^3-1)>1$ (e.g.~$p=7$, $q=19$), then a
cyclic extension of degree $n$ admits only cyclic Hopf-Galois
structures \cite[Theorems 24, 25]{NYJM}. In both cases, non-abelian groups of
order $n$ exist, but a partial analysis of their holomorphs shows
that they cannot arise as the type of a Hopf-Galois structure on a
cyclic extension.
\end{remark}

\section{Abelian Hopf-Galois structures on abelian extensions}

In this final section we describe an alternative approach to Theorem
\ref{nilp-thm} in the case that $\Gamma$ and $G$ are both abelian
(restated as Theorem \ref{ab-case} below). This avoids the use of
Proposition \ref{cent-prop}, and instead is based upon a
result of Caranti, Dalla Volta and Sala \cite{CDVS} which underlies Lemma
\ref{FCC-lemma}. It therefore shows how the ideas in \cite{FCC} extend
to a finite abelian group $\Gamma$ which is not of prime-power order. 

An important ingredient in the proof of Lemma \ref{FCC-lemma} (though
not of the original weaker version in Featherstonhaugh's thesis \cite{F}) is a
correspondence between regular subgroups of $\Hol(G)$ for an abelian
group $G$ and certain multiplication operations $\cdot$ on $G$.
This correspondence was first observed in \cite[Theorem
  1]{CDVS} for vector spaces over a field $F$. The case
$F=\F_p$ (the field of $p$ elements) covers elementary
abelian $p$-groups $G$. It was noted in \cite{FCC} that the same
argument works for any finite $p$-group; indeed, this is what is
required to prove Lemma \ref{FCC-lemma}. It is easily verified that the
argument of \cite{CDVS} is still valid for arbitrary abelian
groups. Here is the result in that setting.  

\begin{lemma}  \label{CDVS-corr}
Let $(G,+)$ be an abelian group with identity element $0$. Then there
is a one-to-one correspondence between regular abelian subgroups $T$
of $\Hol(G)$ and binary operations $\cdot$ on $G$ which make
$(G,+,\cdot)$ into a commutative, associative (non-unital) ring with
the property that every element of $G$ has an inverse under the circle
operation $x \circ y = x + y + x \cdot y$ (so $(G,\circ)$ is an
abelian group, whose identity element is again $0$). Under this
correspondence, the subgroup $T$ of $\Hol(G)$ corresponding to $\cdot$
is $\{ \tau_g \; : \; g\in G\}$, where $\tau_g(x)=g \circ x$ for all
$x \in G$.
\end{lemma}

We next investigate the Sylow subgroups of (the additive group of)
such a ring.

\begin{proposition} \label{fin-ring}
Let $(R,+, \cdot)$ be a finite associative non-unital ring, and for each prime
$p$ dividing its order, let $R_p$ be the Sylow $p$-subgroup of
$(R,+)$. Then $R_p$ is an ideal (and hence a subring) of $R$, and
$R$ is the direct product of its subrings
$R_p$. Moreover, every element of $R$ has an inverse under
$\circ$ if and only if the same is true in each $R_p$.
\end{proposition}
\begin{pf}
Let $r \in R_p$, and let $s \in R$ be arbitrary. If $p^e$ is the
exponent of $R_p$ then, by associativity, $p^e(r \cdot s)=(p^e r)\cdot s=0
\cdot s =0$, so that $r \cdot s \in R_p$. Similarly $s \cdot r \in
R_p$. In particular, if $r \in R_p$ and $s \in R_p$ then $r \cdot s
\in R_p$, and if $r \in R_p$ 
and $s \in R_q$ with $p \neq q$ then $r \cdot s \in R_p \cap R_q$ so
$r \cdot s = 0$. Hence  
$R_p$ is both an ideal and a subring of $R$, and $R$ is the direct
product of its subrings $R_p$.   
Suppose now that every $r \in R$ has a $\circ$-inverse. If $r \in
R_p$ has $\circ$-inverse $s$ in $R$ then 
$s=-r-r\cdot s \in R_p$, so $r$ has $\circ$-inverse $s$ in
$R_p$. Conversely, suppose that $\circ$-inverses exist in each $R_p$. Let $r\in R$. We can
write $r= \sum_p r_p$ with $r_p \in R_p$ for each $p$. If $s_p$ is the
$\circ$-inverse of $r_p$ in $R_p$ then $s= \sum_p s_p$ is the
$\circ$-inverse of $r$ in $R$.  
\end{pf}

\begin{corollary} \label{T-syl}
In Lemma \ref{CDVS-corr}, the Sylow $p$-subgroup $T_p$ of $T$ is $\{
\tau_g \; : \; g \in G_p\}$.
\end{corollary}
\begin{pf}
If $g$, $h \in G_p$ then $g \circ h =g+h +g \cdot h\in G_p$ by Proposition
\ref{fin-ring}. But $\tau_g (\tau_h(x))=g \circ (h \circ x)=(g \circ
h) \circ x= \tau_{g \circ h}(x)$. The non-empty subset $\{\tau_g \; : \; g
\in G_p\}$ of the finite abelian group $T$ is therefore closed under composition, and
hence is a subgroup. Since its cardinality is $|G_p|$ and
$|G|=|T|$, it is the Sylow $p$-subgroup $T_p$.
\end{pf}

\begin{theorem} \label{ab-case}
Let $\Gamma$ and $G$ be abelian groups of order $n$. Then
$$ e(\Gamma,G)= \prod_{p | n} e(\Gamma_p,G_p).  $$
\end{theorem}
\begin{pf}
Let $\beta \colon \Gamma \lra \Hol(G)$ be a regular embedding. Then
$T=\beta(\Gamma) \cong \Gamma$ is a regular subgroup of $\Hol(G)$
which by Lemma \ref{CDVS-corr} gives a multiplication $\cdot$ on $G$
making $G$ into a ring.  Then $T = \{ \tau_g \; : g \in G\}$, where
the $\tau_g$ are defined using the $\circ$-operation obtained from
$\cdot$.  By Proposition \ref{fin-ring}, $G$ is the direct product of
its subrings $G_p$. Since $\circ$-inverses exist in $G$, they exist in
$G_p$, so that the multiplication on $G_p$ corresponds via Lemma
\ref{CDVS-corr} to a regular subgroup $T'_p$ of $\Hol(G_p)$. Writing
elements of $G= \prod_p G_p$ as tuples $g=(g_p)_p$ with $g_p \in G_p$,
we have
$$ \tau_g(x) = g + x + g \cdot x = (g_p + x_p + g_p \cdot x_p )_p $$
for any $x=(x_p)_p \in G$. It follows that $T'_p$ consists of the
restrictions to $G_p$ of the $\tau_{g_p}$ for $g_p \in G_p$. By
Corollary \ref{T-syl}, the $\tau_{g_p}$ are precisely the elements of
the Sylow $p$-subgroup $T_p=\beta(\Gamma_p)$ of $T$. Thus $\beta$
induces a regular embedding $\beta_p \colon \Gamma_p \lra \Hol(G_p)$
for each $p$, where $\beta_p(h)$ for $h \in G_p$ is merely the restriction of
$\beta(h)$ to $G_p$. If we form the product $\beta^* = \left(
\prod_p \beta_p \right) \colon \Gamma \lra \Hol(G)$ as in
(\ref{dp-beta}), then $T^*=\beta^*(\Gamma)$ is a regular subgroup of
$\Hol(G)$ which induces the operation $\cdot$ on each $G_p$. By Lemma
\ref{CDVS-corr} and Proposition \ref{fin-ring} we then have $T^*=T$
and so $\beta^*=\beta$. Thus every regular embedding $\beta$ comes
from a family of regular embeddings $\beta_p$. As in the proof of
Theorem \ref{nilp-thm}, it follows that $e(\Gamma,G)=\prod_p
e(\Gamma_p,G_p)$.  
\end{pf}

\end{document}